\documentclass{amsart}
\usepackage{graphicx}
\usepackage[all]{xy}
\usepackage{amsmath,amssymb,amscd,amsfonts,mathrsfs}
\usepackage[mathcal]{euscript}

\newtheorem{thm}{Theorem}[section]
\newtheorem{conj}[thm]{Conjecture}

\newtheorem{lem}[thm]{Lemma}

\theoremstyle{definition}

\newtheorem{ex}[thm]{Example}
\theoremstyle{remark}
\newtheorem{rem}[thm]{Remark}
\numberwithin{equation}{section}
\newcommand{\R}{\mathbb R}
\newcommand{\eps}{\varepsilon}

\newcommand{\X}{X_0(N)}
\newcommand{\SL}{ SL_2(\Z)}
\newcommand{\Q}{\mathbb{Q}}
\newcommand{\Z}{\mathbb{Z}}
\newcommand{\C}{\mathbb{C}}
\newcommand{\G}{\Gamma_0(N)}
\newcommand{\Gs}{\Gamma^\ast_0(N)}
\newcommand{\g}{\gamma}

\newcommand{\h}{\mathfrak{h}}
\newcommand{\pat}{\partial}
\newcommand{\Sk}{S_{2k}(N)}
\newcommand{\heg}{\mathcal{H}_N^D}
\newcommand{\hegN}{\mathcal{H}_N}

\newcommand{\kd}{\kappa_D}

\newcommand{\inti}{\int\limits_{\I\infty}^{\g(\I\infty)} f(z)}

\newcommand{\I}{i}

\def\Rea{\operatorname{Re}}

\newcommand{\om}{\omega}

\newcommand{\hforms}{\mathcal{Q}_N^D(r)}

\begin{document}

\title[Higher Weight Heegner Points ]{Higher Weight Heegner Points }%
\author{Kimberly Hopkins }%
\address{UT Austin, Department of Mathematics C1200, Austin, TX 78712 }%
\email{khopkins@math.utexas.edu }%

%\date{3/30/09}% remove later!
%\thanks{ }%
\subjclass{ }%
\keywords{modular forms, elliptic curves, Gross-Kohnen-Zagier formula, Shimura-Kohnen correspondence}%

%\date{}%
%\dedicatory{}%
%\commby{}%
% ----------------------------------------------------------------
\begin{abstract}
In this paper we formulate a conjecture which  partially generalizes the Gross-Kohnen-Zagier theorem to higher weight modular forms. For $f\in \Sk$ satisfying certain conditions, we construct a map from the Heegner points of level $N$ to a complex torus, $\C/L_f$, defined by $f$. We define higher weight analogues of Heegner divisors on $\C/L_f$. We conjecture they all lie on a line, and their positions are given by the coefficients of a certain Jacobi form corresponding to $f$. In weight $2$, our map is the modular parametrization map (restricted to Heegner points), and our conjectures are implied by Gross-Kohnen-Zagier. For any weight, we expect that our  map is the Abel-Jacobi map on a certain modular variety, and so our conjectures are consistent with the conjectures of Beilinson-Bloch. We have verified our map is the Abel-Jacobi for weight $4$. We provide numerical evidence to support our conjecture for a variety of examples.
\end{abstract}
\maketitle
% ----------------------------------------------------------------
\section{Introduction }\label{S:intro}

For integers $N, k\geq 1$, let $S_{2k}(N)$ denote the cusp forms of weight $2k$ on the congruence group $\G$. Let $\X$ be the usual modular curve and $J_0(N)$ its Jacobian.  By $D$ we will always mean a negative fundamental discriminant which is a square modulo $4N$. For each   $D$, one can construct a Heegner divisor $y_D$ in $J_0(N)$ and defined over $\Q$. Suppose $f\in S_2(N)$ is any normalized newform whose sign in the functional equation of $L(f,s)$ is $-1$. Then the celebrated theorem of Gross, Kohnen, and Zagier \cite[Theorem C]{GZK} says that, as $D$ varies, the $f$-eigencomponents of the Heegner divisors $y_D$ all `lie on a line\footnote{We will say a subset $X$ of an abelian group $J$ lies on a line if $X \subseteq \Z\cdot x_0$ for some $x_0 \in J$.}'  in the quotient $J_0(N)_f$. Furthermore it says their positions on this line are given by the coefficients of a certain Jacobi form. In particular when $N$ is prime, the positions are the coefficients of a half-integer weight modular form in Shimura correspondence with $f$.

Now suppose $f\in\Sk$ is a normalized newform of weight $2k$ and level $N$. In addition, assume the coefficients in its Fourier series are rational, and the sign in the functional equation of $L(f,s)$ is $-1$. Let $\hegN/\G \subset \X$ denote  the  Heegner points of level $N$. In this paper we construct a map,
\[
    \alpha: \hegN/\G \rightarrow \C/L_f,
\]
where $\C/L_f$ is a complex torus defined by the periods of $f$. Let $h(D)$ denote the class number of the imaginary quadratic field of discriminant $D$. For each  $D$ and fixed choice of its square root $(\bmod2N)$, we get precisely $h(D)$ distinct representatives $\tau_1, \dots, \tau_{h(D)}$ of $\hegN/\G$. Define $(\mathcal{Y}_D)_f = \alpha(\tau_1) + \cdots + \alpha(\tau_{h(D)})$  and define $(y_D)_f = (\mathcal{Y}_D)_f + \overline{(\mathcal{Y}_D)_f}$ in $\C/L_f$. When $k=1$, $\alpha$ is the usual modular parametrization map restricted to Heegner points, and $(y_D)_f$ is equal to the $f$-eigencomponent of $y_D$ in $J_0(N)$ as described in the first paragraph. For $k\geq 1$ we formulate conjectures similar to Gross-Kohnen-Zagier. We predict the $(y_D)_f$ all lie  in a line in $\C/L_f$, that is, there exists a point $y_f \in \C/L_f$ such that
\[
    (y_D)_f = m_D y_f,
\]
up to torsion,  with $m_D \in \Z$. Furthermore we predict the positions $m_D$ on the line are coefficients of a certain Jacobi form corresponding to $f$. In the case when $N$ is prime and $k$ is odd, the $m_D$ should be the coefficients of a weight $(k+1/2)$ modular form in Shimura correspondence with $f$.

We expect our map is equivalent to the Abel Jacobi map  on Kuga-Sato varieties in the following sense. Let $Y=Y^k$ be the Kuga-Sato variety associated to weight $2k$ forms on $\G$. (See \cite[p.$117$]{Zh} for details.)  This is a smooth projective variety over $\Q$ of dimension $2k-1$. Set $\mathcal{Z}^k(Y)_\text{hom}$ to be the nullhomologous codimension $k$ algebraic cycles, and $\text{CH}^k(Y)_{\text{hom}}$ the group of $\mathcal{Z}^k(Y)_\text{hom}$ modulo rational equivalence. Let $\Phi^k$ be the usual $k$-th Abel-Jacobi map,
\[
    \Phi^k: \text{CH}^k(Y)_{\text{hom}} \rightarrow J^k(Y),
\]
where $J^k(Y)$ is the $k$-th intermediate Jacobian of $Y$. Given any normalized newform $f = \sum_{n\geq 1} a_n q^n \in \Sk$ with rational coefficients, there exists an $f$-isotypical component $J_f^k(Y)$ of $J^k(Y)$, and thus an induced map,
\[
\xymatrix{
 \text{CH}^k(Y)_{\text{hom}} \ar[r]^{\Phi^k} \ar@{-->}[rd]_{\Phi^k_f} & J^k(Y) \ar[d]\\
&J_f^k(Y)}
\]
Our expectation is that the image of $\Phi_f^k$ on classes of CM cycles in $\text{CH}^k(Y)_{\text{hom}}$ is equal (up to a constant) to the image of our map $\alpha$ on Heegner points in $\X$. If we assume this is the case, then our conjectures are consistent with the conjectures of Beilinson and Bloch. In this setting they predict
\[
    \text{rank}_\Z \: \text{CH}^k(Y_F)_{\text{hom}} = \text{ord}_{s=k} L_F(H^{2k-1}(Y),s).
 \]
If we assume $\text{ord}_{s=k}L(f,s)=1$, then a refinement of their conjecture predicts the image of $\Phi_f^k$ on CM divisors in $Y_\Q$ should have rank at most $1$ in $J_f^k(Y)$.

We have verified the equivalence of $\alpha$  and $\Phi^2_f$ in the case of weight $4$. For this we used an explicit description of $\Phi^2_f$ on CM cycles given by Schoen in \cite{Sc1}. In fact, in \cite{Sc2} Schoen uses this map to investigate a consequence of Beilinson-Bloch similar to the one described above. For a specific $Y=Y^4$ and $f$ he computes $\Phi_f$ on certain CM divisors in $Y$ defined over the quadratic number field $\Q(i)$. From this he finds numerical evidence that the images lie on a line and their positions are given by a certain weight $5/2$ form corresponding to $f$.

The sections of this paper are divided as follows. In Section \ref{S:mymap} we describe our map and its lattice of periods. In Section \ref{S:conjectures} we give explicit statements of our conjectures. In Section \ref{S:algorithm} we describe the algorithm we created to numerically verify the conjectures in a variety of examples. Note our algorithm could be applied to compute coefficients of half-integer weight modular forms. In sections \ref{S:examples} and \ref{S:moreexamples} we compute some examples and use them to verify our conjectures in two different ways.

%-----------------------------------------------------------------
\section{Higher Weight Heegner Points}\label{S:mymap}

Let  $\h$ denote the upper half-plane. Suppose $f$ is a normalized newform in $\Sk$ having a  Fourier expansion of the form,
\[
    f(\tau) = \sum_{n=1}^\infty a_n q^n, \qquad q = \exp(2\pi i \tau),\: \tau \in \h
\]
with $a_n \in \Q$.

Recall the $L$-function of $f$ is defined by the Dirichlet series,
\[
    L(f,s) = \sum_{n=1}^\infty \frac{a_n}{n^s}, \qquad \Rea(s)>k+1/2,
\]
and has an analytic continuation to all of $\C$. Moreover the function $\Lambda(f,s) = N^{s/2} (2\pi)^{-s} \Gamma(s) L(f,s)$  satisfies the functional equation,
\[
    \Lambda(f,s) = \eps \Lambda(f, 2k-s),
\]
where $\eps = \pm 1$ is   the sign of the functional equation of $L(f,s)$.

For each prime divisor $p$ of $N$, let $q = p^\ell$, $\ell \in \mathbb{N}$ such that $\gcd(q,N/q)=1$ and set $\om_q = \begin{smat}
{qx_0}{1}{Ny_0}{q}
\end{smat}$, for some $x_0,y_0\in \Z$, with $qx_0 - (N/q)y_0 =1$ . Define $\Gs$ to be the group generated by $\G$ and each $\om_q$. Let $S$ be a set of generators for $\Gs$. Define the period integrals of $f$ for the set $S$ by,
\[
    \mathcal{P} = \bigg\{ (2\pi i)^{k} \inti z^m dz \::\: m\in \{0,\dots, 2k-2\},\: \g\in S \bigg\}\subseteq \C.
\]
These are sometimes referred to as Shimura integrals. It is straightforward to see that every  integral of the form,
\[
    (2\pi i)^{k} \inti z^m dz, \qquad \g \in \Gs,\: 0\leq m \leq 2k-2
\]
is in an integral linear combination of elements in $\mathcal{P}$. (See \cite[Section $8.2$]{Sh}, for example). In fact, the $\Z$-module generated by $\mathcal{P}$ forms a lattice,
\begin{lem}
$L:=\text{Span}_\Z(\mathcal{P})$ is a lattice in $\C$.
\end{lem}
\begin{proof}

By theorems of Razar \cite[Theorem $4$]{Ra} and \v{S}okurov \cite[Lemma $5.6$]{So}, the set $\mathcal{P} $ is contained in some lattice. Hence $L$ is of rank $\leq 2$. To show its rank is $2$, it suffices to show there exist  nonzero complex numbers $u^+, u^- \in L$ with $u^+ \in \R$ and $u^-\in i\R$.

Suppose $m$ is a prime not dividing $N$, and  $\chi$  a primitive Dirichlet character modulo $m$.  Define $(f\otimes \chi) := \sum_{n\geq 1} \chi(n) a_n q^n$, and $L(f\otimes \chi, s)$ to be its Dirichlet series. Let $\Lambda(f\otimes \chi,s) = (2\pi)^{-s} (Nm^2)^{s/2} \Gamma(s) L(f\otimes \chi, s)$. Then for $\Rea(s)>k+1/2$, we have
\begin{equation}\label{E:lamfchi}
    i^s (Nm^2)^{-s/2} \Lambda(f\otimes \chi, s) =  \int\limits_0^{i\infty} (f\otimes\chi) (z) z^s \frac{dz}{z}.
\end{equation}

Let $g(\chi)$ denote the Gauss sum associated to $\chi$. Then an expression for $\chi$ in terms of the additive characters is given by,
\[
    \chi(n) = m^{-1} g(\chi) \sum_{u\bmod m} \bar{\chi}(-u) e^{2\pi i nu/m}.
\]
So
\[
    (f\otimes \chi)(\tau) = m^{-1} g(\chi) \sum_{u\bmod m} \bar{\chi}(-u) f(z + u/m).
\]
Substituting this into \eqref{E:lamfchi} gives
\[
    i^s (Nm^2)^{-s/2} \Lambda(f\otimes \chi, s) = m^{-1} g(\chi) \sum_{u\bmod m} \bar{\chi}(-u)  \int\limits_0^{i\infty} f(z+u/m) z^s \frac{dz}{z},
\]
and replacing $z$ by $z-u/m$ and rearranging implies
\[
    i^{-s} g(\chi)^{-1} N^{-s/2}  \Lambda(f\otimes \chi, s) = (-1)^{s-1}\sum_{u\bmod m} \bar{\chi}(-u) \int\limits_{i\infty}^{u/m} f(z) (mz-u)^{s-1} dz.
\]

Now let $s=2k-1$ in the above equation, and multiply both sides by $(2\pi i)^k$. In addition suppose $\chi$ is a  quadratic Dirichlet character modulo $m$. If $m\equiv 3\bmod 4$, then $g(\chi) = i\sqrt{m}$, and if $m\equiv 1\bmod 4$ then $g(\chi) = \sqrt{m}$. Hence since $\Lambda(f\otimes \chi,2k-1)$ is real-valued and nonzero, the right hand side of this equation is either purely real or purely imaginary depending on the choice of $m$. Then this proves the lemma since the right hand side is in $L$ for any $m$.

\end{proof}

Let $D<0$ be a fundamental discriminant, and assume $D$ is  a square modulo $4N$.  Fix a residue class $r\bmod 2N$ satisfying $D\equiv r^2\bmod 4N$. Then
\[
    \hforms := \big\{ [A,B,C] :  A>0, B, C \in \Z, D= B^2-4AC,  A \equiv 0\bmod N, B\equiv r\bmod 2N \big\}.
\]
corresponds to a subset of the positive definite binary quadratic forms of discriminant $D$.  We define $\heg(r)$  to be the roots in $\h$ of $\hforms$,
\[
    \heg(r) := \bigg\{ \tau = \frac{-B+\sqrt{D}}{2A} : [A,B, C] \in \hforms, C= \frac{|D|+B^2}{4A} \bigg\}.
\]
$\G$ preserves $\heg(r)$, and the classes of $\heg(r)/\G$ are in bijection with the classes of reduced binary quadratic forms of discriminant $D$.  We will call $\heg(r)/\G$ the set of Heegner points of level $N$, discriminant $D$, and root $r$. Define  $\hegN$ to be the union of $\heg(r)$ over all $D, r$, and so $\hegN/\G$ are the Heegner points of level $N$.

For each $\tau = \frac{-B+\sqrt{D}}{2A} \in \heg(r)$, set $Q_\tau(z) := Az^2 +Bz+C$. We now define a function $\alpha =\alpha_f: \hegN \rightarrow \C$ by
\[
    \alpha(\tau) :=  (2\pi i )^{k} \int_{i\infty}^\tau f(z) Q_\tau(z)^{k-1}dz.
\]

\begin{lem}\label{L:welldef}
The map $\alpha$ induces a well-defined map (which we will also denote by $\alpha$),
\[
    \alpha: \hegN/\G \rightarrow \C/L.
\]
\end{lem}
\begin{proof}
For any $\tau \in \hegN$ of discriminant $D$ and $\g\in \G$, we will show
\[
    \alpha(\g \tau) - \alpha(\tau) = \:(2\pi i)^{k}\cdot\!\!\!\inti  Q_{\g \tau}(z)^{k-1} dz.
\]
Since $Q_{\g \tau}(z)$ has integer coefficients, this will imply $\alpha(\g \tau) - \alpha(\tau) \in L$ for all $\g \in \G$.

Let $\g = \begin{smat}
{a}{b}{c}{d}
\end{smat} \in \G$. Then
\begin{align*}
 & \alpha(\g \tau) - (2\pi i)^{k}\cdot\!\!\! \inti \:  Q_{\g \tau}(z)^{k-1}  dz \\
&= (2\pi i)^{k}\cdot\!\!\! \int\limits_{\g(i\infty)}^{\g\tau} f( z)  Q_{\g \tau}(z)^{k-1}  dz\\
            &= (2\pi i)^{k}\cdot\!\!\! \int\limits_{i\infty}^{\tau} f(\g z) Q_{\g \tau}(\g z)^{k-1} d(\g z)\\
            &= \alpha(\tau),
\end{align*}
where in the last equality we used $f(\g z) = (cz+d)^{2k} f(z)$, $Q_{\g \tau}(z) = (-cz+a)^2 Q_\tau(z)$, and $d(\g z) = (cz+d)^{-2} dz$.

\end{proof}

%---------------------------------------------------------------

\section{Conjectures}\label{S:conjectures}

Let $\{\tau_1, \dots, \tau_{h(D)}\} \in \heg(r)$ be any set of  distinct class representatives of $\heg(r)/\G$. Define
\[
    P_{D,r} := \sum_{i=1}^{h(D)} \tau_i  \in \text{Div}(\X),
\]
where $\text{Div}(\X)$ denotes the group of divisors on $\X$.  If $D=-3$ (resp. $D=-4$ ), scale $P_{D,r}$ by  $1/3$ (resp. $1/2$). Extend $\alpha$ to $P_D$ by linearity and define
\[
   (y_{D,r})_f =  \alpha(P_{D,r}) + \overline{\alpha(P_{D,r})} \in \C/L.
\]
Here, bar denotes complex conjugation in $\C$.  We write $y_{D,r}$ or $y_D$ for $(y_{D,r})_f$, and $P_D$ for $P_{D,r}$ when the context of $f$, $r$ is clear.

By the actions of complex conjugation and Atkin-Lehner on $\hegN$, we have
\[
\overline{\alpha(P_{D,r})} = -\eps \alpha(P_{D,r}),
\]
where  $\eps$ is the sign of the functional equation of $L(f,s)$. Thus if $\eps=+1$, then $y_{D,r}$ are in $L$ for all $D,r$. This is, in some sense, the trivial case. Hence we restrict our attention to the case when $\eps=-1$.

Conjectures \ref{C:conj1} and \ref{C:conj2} give a partial generalization of the Gross-Kohnen-Zagier theorem to higher weights.

\begin{conj}\label{C:conj1}
Let $f = \sum_{n\geq 1} a_n q^n \in \Sk$ be a normalized newform with rational coefficients, and assume $\eps=-1$ and $L'(f,k) \neq 0$. Then for all fundamental discriminants $D<0$ and $r\bmod 2N$ with $D\equiv r^2 \bmod 4N$, there exist integers $m_{D,r}$ such that
\[
    ty_{D,r} = m_{D,r} y_f \qquad \text{in}\: \: \C/L,
\]
where $y_f \in \C/L$ and $t\in \Z$ are both nonzero and independent of $D$ and $r$.

\end{conj}

\begin{rem}
Equivalently we could say $y_{D,r} = m_{D,r} y_f$ up to a $t$-torsion element in $\C/L$.
\end{rem}

To state the second conjecture we will need to use Jacobi forms. (See \cite{Ei} for background). Let $J_{2k,N}$ denote the set of all  Jacobi forms of weight $2k$ and index $N$. Then such a $\phi \in J_{2k,N}$ is a function $\phi: \h \times \C \rightarrow \C$, which satisfies the transformation law
\[
    \phi \bigg( \frac{a\tau+b}{c\tau +d}, \frac{z}{c\tau +d} \bigg) = (c\tau+d)^{2k} e^{2\pi i N \frac{cz^2}{c\tau+d}} \phi(\tau, z),
\]
for all $\begin{smat}
{a}{b}{c}{d}
\end{smat} \in \SL$, and has a Fourier expansion of the form
\begin{equation}\label{E:jacfour}
    \phi(\tau, z) = \sum_{\begin{subarray}{1} n, r \in \Z\\
       r^2 \leq 4Nn \end{subarray} } c(n,r) q^n \zeta^r, \qquad q= e^{2\pi i \tau}, \:\zeta = e^{2\pi i z}.
\end{equation}
The coefficient $c(n,r)$ depends only on $r^2 - 4Nn$ and on the class $r\bmod 2N$.

Suppose $f \in \Sk$ is a normalized newform with $\eps=-1$. Then by \cite{SkZa}, there exists a non-zero Jacobi cusp form $\phi_f \in J_{k+1, N}$ which is unique up to scalar multiple and has the same eigenvalues as $f$ under the Hecke operators $T_m$ for $m, N$ coprime. We predict that the coefficients of $\phi_f$ are related to the $m_D$ from above in the following way,
\begin{conj}\label{C:conj2}
Let $f = \sum_{n\geq 1} a_n q^n \in \Sk$ be a normalized newform with rational coefficients, and assume $\eps=-1$ and $L'(f,k) \neq 0$. Assume Conjecture \ref{C:conj1}. Then
\[
    m_{D,r} = c(n,r)
\]
where $n = \frac{|D|+r^2}{4N}$ and  $c(n,r)$ is the $(n,r)$-th coefficient of the Jacobi form $  \phi_f \in J_{k+1,N}$.
\end{conj}

\begin{rem}
When $k=2$, the points $(y_{D,r})_f$ and $y_f$ are the same as those defined in \cite{GZK}, and both of our conjectures are implied by Theorem C of their paper. (Actually their theorem is only for $D$ coprime to $2N$ but they say the result remains `doubtless true' with this hypothesis removed. See \cite{Ha} and \cite{Bo} for more details.) Particular to weight $2$ is the fact that $\C/L$ is defined over $\Q$ and so $y_D$ is a rational point on the elliptic curve $E_f \simeq \C/L$.  In contrast, we should stress that for weight $k>2$, the elliptic curve $E \simeq \C/L$ is not expected to be defined over any number field. For instance, the $j$-invariants for our examples all appear to be transcendental over $\Q$.
\end{rem}

\begin{rem} For $N=1$ or a prime, and $k$ odd we can state Conjecture \ref{C:conj2} in terms of modular forms of half-integer weight. Specifically, let $\phi\in J_{k+1}(N)$ be a Jacobi form with a  Fourier expansion as in \eqref{E:jacfour}, and set
\[ g(\tau) = \sum_{M=0}^\infty c(M) q^M, \qquad q=e^{2\pi i \tau} \] where $c(M)$ is defined by,
\[ c(M) := \left\{ \begin{array}{ll} c\bigg( \frac{M+r^2}{4N}, r\bigg) & \hbox{if $M\equiv -r^2 \bmod 4N$ for any $r\in \Z$;} \\ 0, & \hbox{otherwise.} \end{array} \right.
 \]
This function is well-defined because $c(n,r)$ depends only on $r^2-4nN$ when $N=1$ or a prime, and $k$ is odd. Then by \cite[p.$69$]{Ei}, $g $ is in $M_{k+1/2}(4N)$, the space of modular forms of weight $k+1/2$ and level $4N$. In addition, if $f\in \Sk$ is a normalized newform with $\eps=-1$, then the form $g$ defined by $\phi_f$ is in Shimura correspondence with $f$.
\end{rem}

\section{Algorithm}\label{S:algorithm}

Let $f =\sum_{n\geq 1} a_nq^n \in \Sk$ be a normalized newform with rational Fourier coefficients.  The sign $\eps$ of the functional equation of $L(f,s)$ can be computed with the  identity,
\[
    f\bigg(\frac{-1}{Nz}\bigg) = (-1)^{k} \eps N^{k} z^{2k} f(z)
\]
given by the action of the Fricke involution of level $N$ on $f$. We will only consider $f$ such that $\eps=-1$ and $L'(f,k)\neq 0$.

The first step is to find a basis of our lattice $L$, which  is the $\Z$-module generated by the periods $\mathcal{P}$ as described above. Suppose $p_1, p_2, p_3$ are three periods in $\mathcal{P}$. Since $L$ has rank $2$, these are linearly dependent over $\Z$, that is
\[
    a_1 p_1 + a_2 p_2 + a_3 p_3=0, \quad\text{for some}\quad a_i \in \Z.
\]
We may assume $\gcd(a_1, a_2, a_3)=1$. Let $d = \gcd(a_1, a_2)$, then there exist integers $x,y \in \Z$ such that $xa_1 + ya_2 =d$. Similarly $\gcd(d,a_3)=1$ so there exist integers $u,v \in \Z$ such that $ud+ va_3 = 1$. Define the matrix $M$ by,
\[
M =  \left(
       \begin{array}{ccc}
         a_1 & a_2 & a_3 \\
         -y & x & 0 \\
         -va_1/d & -va_2/d & u \\
       \end{array}
     \right).
\]
Observe $M \in GL_3(\Z)$ and $M\cdot {}^T(p_1,p_2,p_3) = {}^T(0, -yp_1+xp_2, -va_1p_1/d - va_2p_2/d + up_3)$. Hence $-yp_1+xp_2$ and $  -va_1p_1/d - va_2p_2/d + up_3$ are a basis for the $\Z$-module generated by $p_1, p_2, p_3$.

We would also like our basis elements to have small norm. Given a basis $\om_1,\om_2$ of a lattice, its norm form is a real bilinear quadratic form defined by the matrix,
\[
    B = \left(
          \begin{array}{cc}
            2|\om_1|^2 & 2\Rea(\om_1 \bar{\om}_2) \\
            2\Rea(\om_1 \bar{\om}_2) & 2|\om_2|^2 \\
          \end{array}
        \right).
\]
Thus it is equivalent to a reduced form of the same discriminant, that is, there exists $U \in \SL$ such that
\[
    {}^T U BU = \left(
                  \begin{array}{cc}
                    2\alpha & \beta \\
                    \beta & 2\gamma \\
                  \end{array}
                \right), \qquad \alpha, \beta, \gamma \in \R,
\]
with $|\beta|\leq \alpha \leq \gamma$ and $\beta\geq 0$ if either $|\beta| = \alpha$ or $\alpha=\gamma$. Hence $(\om_1', \om_2') := (\om_1,\om_2)U$ is a `reduced' basis. For a basis of all of $L$ we simply apply this process iteratively on the elements of $\mathcal{P}$.

In fact it is not hard to see that $L$ is a real lattice, that is, $\bar{L} = L$. Thus given a basis $\om_1, \om_2$ of $L$, we may assume $\om_1\in i\R$, and therefore $\tau:= \om_2/\om_1$ has real part equal to either $0$ or $1/2$. This implies $\Rea(L) = \Rea(\om_2)$ which will help simplify our computations.

To actually compute the elements in $\mathcal{P}$ we need to split the path from $(i\infty)$ to $\g(i\infty)$ of integration at some point $\tau \in \h$ which gives,
\[
    \inti z^m dz  = \int\limits_{\I\infty}^{\g(\tau)}  f(z) z^m dz - \int\limits_{\I\infty}^{\tau}  f(z) (az+b)^m (cz+d)^{2k-2-m} dz,
\]
for $\g = \begin{smat}
{a}{b}{c}{d}
\end{smat} \in\G$. We choose $\tau$ to be a point at which $f$ has good convergence. To compute integrals of the form,
\[
    \int\limits_{\I\infty}^{\tau} f(z) z^m dz,
\]
we use repeated integration by parts to get the formula
\begin{equation}\label{E:compzpwr}
 \int\limits_{\I\infty}^{\tau} \!\! f(z) z^m dz = m! \,(-1)^m \sum_{j=-1}^{m-1} \frac{(-1)^{j+1}}{(j+1)!} \tau^{j+1} f_{m-j}(\tau),
\end{equation}
where $f_\ell(\tau)$ is defined to be the $\ell$-fold integral of $f$ evaluated at $\tau \in \h$, that is,
\[
    f_\ell(\tau) = \frac{1}{(2\pi i)^\ell} \sum_{n\geq1} \frac{a_n}{n^\ell} q^n, \qquad q = \exp(2\pi i \tau)
\]
which is well-defined for any $0\leq \ell \leq 2k-1$.

The next task is to compute $\alpha(\tau)$ for $\tau \in \hegN$. We could do this using \eqref{E:compzpwr}, but it is computationally faster to use the following identity for $\alpha$. Recall the modular differential operator,
\[
    \pat_m := \frac{1}{2\pi i} \frac{d}{dz} - \frac{m}{4\pi y},\qquad z=x+iy\in \h,
\]
for any integer $m$. Define $\pat_{m}^\ell(f) := \pat_{m+2(\ell-1)}\circ\cdots\circ \pat_{m+2}\circ\pat_m(f)$ to be the composition of the $\ell$ operators $\pat_m, \pat_{m+2}, \dots, \pat_{m+2(\ell-1)}$. Then a straightforward combinatorial argument  yields the following identity, whose proof we will omit,
\begin{lem}\label{L:altalph}
Let $\tau$ be a Heegner point of level $N$ and discriminant $D$. Then
\[
    \alpha(\tau)  = \kd \cdot\pat_{-2k+2}^{k-1}\circ f_{2k-1}(\tau),
\]
where $\kd =(k-1)! \, (2\pi i)^{k} (2\pi \sqrt{|D|} )^{k-1}$ is a constant depending only on $D$ and $2k$.
\end{lem}
A closed formula for $\pat_m^\ell$ (see \cite{Vi} for example) allows us to write $\alpha$ as
\begin{equation}\label{E:compal}
    \alpha(\tau) = \kd(2\pi i) \bigg( \frac{-y}{\pi}\bigg)^{k} \sum_{n\geq 1}  p\bigg(\frac{2k}{2}, \frac{1}{4\pi y n}  \bigg) a_n q^n,
\end{equation}
where $p(m,x)$, is the polynomial,
\[
    p(m,x) = \sum_{\ell=m}^{2m-1} \binom{m-1}{2m-1-\ell}\frac{(\ell -1 )!}{(m-1)!} x^{\ell}, \qquad m\in \Z,\: x\in \R.
\]
 We compute $\alpha(\tau)$ using \eqref{E:compal}. Also notice that Lemma \ref{L:altalph} perhaps provides further insight into why the map $ \hegN \rightarrow \C/L$ inducing $\alpha$ is invariant under $\G$. Loosely speaking, this is because integrating $f$ $(2k-1)$-times lowers its weight by $2(2k-1)$ and $\pat_{-2k+2}^{k-1}$ increases its weight by $2(k-1)$ to get something morally of weight $0$.

Given a set of Heegner point representatives of level $N$, discriminant $D$, and root $r$, we can use the above to compute $y_{D,r}$. Verifying the first conjecture for each $D, r$ then amounts to choosing a complex number $y_f$, and an integer $t$, both non-zero, and showing the linear dependence,
\begin{equation}\label{E:testeq}
    \Rea(y_{D,r})- m_{D,r} \Rea(y_f) + n_{D,r} \Rea(\om_2)/t = 0
\end{equation}
for some integers $m_{D,r}, n_{D,r}$. The second conjecture consists of comparing the coefficients $m_{D,r}$ of $y_f$ we get above with the Jacobi form coefficients of the form $\phi_f$.

%----------------------------------------------------------------------

\section{Examples}\label{S:examples}

The Fourier coefficients of the forms in these examples were computed using SAGE \cite{sage}. The rest of the calculations were done in PARI/GP \cite{PARI2}.

We will always take a set of generators for $\G$ which includes the translation matrix $T = \begin{smat}
{1}{1}{0}{1}
\end{smat}$ but no other matrix whose $(2,1)$ entry is $0$. The period integrals for $T$ are always $0$ since $i\infty$ is its fixed fixed point, hence we can exclude it from our computations of $\mathcal{P}$. In addition the $(2\pi)^k$ factor in the definitions of $y_D$ and $L$ is left off from the computations, since it is just a scaling factor and requires unnecessary extra precision.

For each example below, we list the number of digits of precision and the number $M$ of terms of $f$ we used. Below that is a set of generators we chose for $\Gamma_0^\ast(N)$ and the bases, $\om_1, \om_2$, we got for $L$ from computing $\mathcal{P}$ and applying the lattice reduction algorithm explained in Section \ref{S:algorithm}. We then provide a table listing the $m_D$ which satisfy equation \eqref{E:testeq} for $t$, $y_f$ of our choosing, and $D$ less than some bound. Without getting into details, the precision we chose depended on the size of the $M$-th term of $f$ and on the a priori knowledge of the size of the coefficients satisfying \eqref{E:testeq}.

\begin{ex}$2k=10$, $N=3$.
The space of cuspidal newforms of weight $10$ and level $3$ has dimension $2$, but only one form has $\eps=-1$. The first few terms of it are
\[
f = q - 36q^2 - 81q^3 + 784q^4 - 1314q^5 + 2916q^6 - 4480q^7 - 9792q^8 + \cdots
\]

\[
\renewcommand{\arraystretch}{1.4}
\begin{array}{ll}
\text{Precision} & 60\\
\text{Number of terms} & 100\\
\Gamma^\ast_0(3) & \big< T, \begin{smat}
{-1}{1}{-3}{2}
\end{smat}, \om_{3} = \begin{smat}
{0}{-1}{3}{0}
\end{smat} \big>\\
\om_1 & -i\cdot0.00088850361439085\dots\\
\om_2 &    0.00002189032158611\dots \\
y_f    &   y_{-8}/2 \\
t      &   1 \\

\end{array}
\]

%Table 1 should appear here.

\begin{table}[h]
\begin{center}
\begin{tabular}{|rr | rr|}\hline
\rule[-6pt]{0pt}{19pt} $|D|$ & $m_D$  &$|D|$ & $m_D$ \\\hline
8 &  2 &  104 &  380   \\
11 &  -5 &   107 &  -507   \\
20 &  8 &   116 &  -40  \\
23 &  8 &   119 &  -560  \\
35 &  42 &   131 &  235   \\
47 &  -48 &   143 &  -376   \\
56 &  0 &   152 &  -364   \\
59 &  -155 &   155 &  -64   \\
68 &  160 &   164 &  -1440   \\
71 &  40 &   167 &  1528   \\
83 &  353 &   179 &  2635   \\
95 &  280 &   191 &  -400   \\\hline
\end{tabular}
\vspace{.1in}
\caption{$f \in S_{10}(3)$. List of $D$, $m_D$ such that $ y_D - m_Dy_f \in L$ for $|D|<200$.}\label{Ta:T1}
\end{center}
\end{table}

The $m_D$ in Table \ref{Ta:T1} give, up to scalar multiple, the coefficients of the weight $11/2$ level $12$ modular form found in \cite[p. $144$]{Ei}. Note we can use the theorems of Waldspurger to get information about the values $L(f,D,k)$ from this table. For example, $L(f,-56,5)=0$.

\end{ex}

%----------------------------------------------------
\begin{ex}$2k=18$, $N=1$.

The weight $18$ level $1$ eigenform in $S_{18}(1)$ has the closed form
\[
    f(z) = \frac{-E_6^3(z) +  E_4^3(z)E_6(z)}{1728},
\]
where $E_{2k}(z)$ is the normalized weight $2k$ Eisenstein series.

\[
\renewcommand{\arraystretch}{1.4}
\begin{array}{ll}
\text{Precision} & 200 \\
\text{Number of terms} & 100 \\
\Gamma^\ast_0(1) = SL_2(\Z) & \big< T,S = \om_1 = \begin{smat}
{0}{-1}{1}{0}
\end{smat} \big>\\
\om_1 &  i\cdot0.001831876775870191761\dots\\
\om_2 &   0.000000000519923858624\dots   \\
y_f    &    y_{-3}/3  \\
t      &   1

\end{array}
\]

%Table 2 should appear here

\begin{table}[h]
\begin{center}
\begin{tabular}{|rr | rr|}\hline
\rule[-6pt]{0pt}{19pt}
$|D|$ & $m_D$ & $|D|$ & $m_D$ \\\hline
3 & 1 &51 & 108102 \\
4 & -2 &52 & -93704 \\
7 & -16 &55 & -22000 \\
8 & 36 &56 & 80784 \\
11 & 99 &59 & -281943 \\
15 & -240 &67 & 659651 \\
19 & -253 &68 & 193392 \\
20 & -1800 &71 & -84816 \\
23 & 2736 &79 & -109088 \\
24 & -1464 &83 & -22455 \\
31 & -6816 &84 & -484368 \\
35 & 27270 &87 & 1050768 \\
39 & -6864 &88 & 143176 \\
40 & 39880 &91 & 195910 \\
43 & -66013 &95 & -370800 \\
47 & 44064 & &  \\\hline

\end{tabular}
\vspace{.1in}
\caption{$f\in S_{18}(1)$. List of $D$, $m_D$ such that $ y_D - m_Dy_f \in L$ for $|D|<100$.}\label{Ta:T2}
\end{center}
\end{table}

The $m_D$ in Table \ref{Ta:T2} are identical to the coefficients of the weight $19/2$ level $4$ half-integer weight form in \cite[p.$141$]{Ei}, which is in Shimura correspondence with $f$.

\end{ex}

%------------------------------------------------------
\begin{ex}$2k=4$, $N=13$

The dimension of the new cuspidal subspace is $3$ in this case, but only one has integer coefficients in its $q$-expansion.
\[
    f = q - 5 q^2 - 7 q^3 + 17 q^4 - 7 q^5 + 35 q^6 - 13 q^7 - 45 q^8 + 22 q^9 + \cdots
\]

\[
\renewcommand{\arraystretch}{1.4}
\begin{array}{ll}
\text{Precision} &  28 \\
\text{Number of terms} & 250  \\
\Gamma^\ast_0(13)  & \big< T, \begin{smat}
{8}{-5}{13}{-8}
\end{smat}, \begin{smat}
{-3}{1}{-13}{4}
\end{smat}, \begin{smat}
{5}{-2}{13}{-5}
\end{smat} ,
\begin{smat}
{-9}{7}{ -13}{10}
\end{smat}, \om_{13} = \begin{smat}
{0}{-1}{13}{0}
\end{smat} \big>\\
\om_1 & i\cdot0.003124357726009878347400865279 \dots \\
\om_2 &  -0.04271662498543992056668379773\dots \\
       & \quad - i\cdot0.001562178863004939178984383052\dots \\
y_f    &   y_{-3}/3  \\
t      &    3

\end{array}
\]

Notice this is the first example of a nonsquare lattice. In fact $\om_2/\om_1 = -0.5000\dots + i\cdot13.67212999\dots$ so $\Rea(\om_2/\om_1)=1/2$ as explained earlier. This is also the first example where the choice of $r$ matters, since $k=2$ is not odd. For each $D$, we chose $r$ in the interval $0<r<13$. In addition this is our only example where $t>1$.

\begin{table}[h]
\begin{center}
\begin{tabular}{|rr | rr|}\hline
\rule[-6pt]{0pt}{19pt}
$|D|$ & $m_{D,r}$ & $|D|$ & $m_{D,r}$ \\\hline
3&1&107&4 \\
4&-1&116&-8 \\
23&2& 120&-13 \\
35&-7&127&14 \\
40&3&131&-3 \\
43&-17&139&29 \\
51&9&152&2 \\
55&-6&155&22 \\
56&1&159&-6 \\
68&-5&168&-21 \\
79&4&179&-17 \\
87&-6&183&-2 \\
88&10&191&-10 \\
95&4&199&4 \\
103&-8 & &\\\hline
\end{tabular}
\vspace{.1in}
\caption{$f\in S_{4}(13)$. List of $D$, $m_{D,r}$ such that $ ty_{D,r} - tm_{D,r}y_f \in L$ with $t=3$, for $|D|<200$.}\label{Ta:T3}
\end{center}
\end{table}

A closed form expression for the weight $3$ index $13$ Jacobi form $\phi = \phi_f$ corresponding to $f$ was provided to us by Nils Skoruppa,
\begin{equation*}
\phi(\tau, z) = \vartheta_1^5 \vartheta_2^3 \vartheta_3 / \eta^3
\end{equation*}
Here $\eta$ is the usual Dedekind eta-function, $\eta = q^{1/24}\prod_{n\geq1} (1-q^n)$ with $q=e^{2\pi i \tau}$, and $\vartheta_a = \sum_{r \in \Z} \left( \frac{-4}r\right)
q^{\frac{r^2 } 8} \zeta^{\frac{ar} 2}$ for $a=1,2,3$, $\zeta = e^{2\pi iz}$. (This has a nice product expansion using Jacobi's triple product identity.)

We verify that the $(n,r)$-th coefficient  $c(n,r)$  in the Fourier expansion of $\phi$ is identically equal to the $m_{D,r}$ in Table \ref{Ta:T3} for $|D|<200$.

\end{ex}

\section{More Examples}\label{S:moreexamples}

The coefficients of Jacobi forms are difficult to compute, in particular for the cases when $N$ is composite or when $k$ is even. We chose the previous examples in part because the Fourier coefficients for their Jacobi forms already existed, thanks to the work of Zagier, Eichler, and Skoruppa mentioned above. However, given any weight and level, we can still provide convincing evidence for our conjecture without knowing the exact coefficients of its Jacobi form. This is done using a refinement of Waldspurger \cite{Wa} given in \cite[p.$527$]{GZK}.

Specifically, let $f\in \Sk$ be a normalized newform with $\eps=-1$. Let $\phi = \phi_f \in J_{k+1,N}$, with Fourier coefficients denoted by $c(n,r)$, be the Jacobi form corresponding to $f$ as described in Section \ref{S:conjectures}. For a fundamental discriminant $D$ with $\gcd(D,N)=1$ and   square root $r$ modulo $4N$, \cite[Corollary $1$]{GZK} says
\[
   |D|^{k-1/2} L(f,D,k) \doteq  |c(n,r)|^2;
\]
 here  $L(f,D,s)$ is $L$-series of $f$ twisted by $D$, and $n\in \Z$ satisfies $D=r^2 -4Nn$. By $\doteq$ we mean equality up to a nonzero factor depending on $N ,2k, f,$ and $\phi$, but independent of $D$. (Gross-Kohnen-Zagier give this constant explicitly in their paper, but for us it is unnecessary.)

Thus given two such discriminants $D_i = r_i^2-4Nn_i$, $i=1,2$, we have
\[
    \frac{ |D_1|^{k-1/2} L(f,D_1,k)}{ |D_2|^{k-1/2} L(f,D_2,k)} = \frac{|c(n_1,r_1)|^2}{|c(n_2,r_2)|^2}.
\]
Hence by computing central values of twisted $L$-functions of $f$, we can test if ratios of squares of our $m_{D_i,r_i}$ are equal to those of $c(n_i,r_i)$.

For the examples below we have the same format as the previous examples along with a fixed choice of discriminant $D_1$ for which we verified explicitly,
\[
    \frac{ |D_1|^{k-1/2} L(f,D_1,k)}{ |D|^{k-1/2} L(f,D,k)} = \frac{m_{D_1,r}^2}{m_{D,r}^2}
\]
for all $D$ coprime to $N$ less than a certain bound.

\begin{ex}$2k=4$, $N=21$.

The dimension of the new cuspidal subspace of $S_4(21)$ is $4$. We chose
\[
    f =  q - 3 q^2 - 3 q^3 + q^4 - 18 q^5 + 9 q^6 + 7 q^7 +\cdots.
\]

\[
\renewcommand{\arraystretch}{1.4}
\begin{array}{ll}
\text{Precision} & 40  \\
\text{Number of terms} &  500 \\
\Gamma^\ast_0(21)  & \big< T, \left(
  \begin{smallmatrix}
    -4 & 1 \\
    -21 & 5 \\
  \end{smallmatrix}
\right),
\left(
  \begin{smallmatrix}
   11 & -5 \\
    42 & -19 \\
  \end{smallmatrix}
\right),
\left(
  \begin{smallmatrix}
    13 & -9 \\
     42 & -29 \\
  \end{smallmatrix}
\right),
\left(
  \begin{smallmatrix}
    8 & -5 \\
    21 & -13 \\
  \end{smallmatrix}
\right),
\left(
  \begin{smallmatrix}
    26 &-19 \\
    63& -46 \\
  \end{smallmatrix}
\right),
\left(
  \begin{smallmatrix}
    -16 & 13 \\
    -21 &  17 \\
  \end{smallmatrix}
\right)   \big>\\
\om_1 &  i\cdot 0.012130626847574141\dots \\
\om_2 &     -0.03257318919429172\dots \\
y_f    &   y_{-3}  \\
t      &   1 \\
D_1   & -20

\end{array}
\]

For a consistent choice of each $r$  we chose the first positive residue modulo $2N$ which satisfies $D\equiv r^2\bmod 4N$ for each $D$.

\begin{table}[h]
\begin{center}
\begin{tabular}{|rr | rr|}\hline
\rule[-6pt]{0pt}{19pt}
$|D|$ & $m_{D,r}$ & $|D|$ & $m_{D,r}$ \\\hline
3 & 1 & 111 & 4 \\
20 & -1 & 119 & 0 \\
24 & -1 &131 & 3 \\
35 & 0 & 132 & 8 \\
47 & 2 & 143 & 2 \\
56 & 0 &152 & -7 \\
59 & 1 &159 & 0 \\
68 & -2 & 164 & -2 \\
83 & 5 & 167 & 4 \\
84 & 0 &168 & 0 \\
87 & -4 & 195 & 8 \\
104 & -3 & &\\\hline
\end{tabular}
\vspace{.1in}
\caption{$f\in S_{4}(21)$. List of $D$, $m_{D,r}$ such that $ y_{D,r} - m_{D,r}y_f \in L$ for $|D|<200$.}\label{Ta:T4}
\end{center}
\end{table}

\end{ex}

%------------------------------------------------------------------------------------------

\begin{ex}$2k=12$, $N=4$.

The space of new cuspforms in $S_{12}(4)$ is spanned by one normalized newform whose Fourier series begins with,
\[
 f = q - 516 q^3 - 10530 q^5 + 49304 q^7 + 89109 q^9 - 309420 q^11 + \cdots.
\]

\[
\renewcommand{\arraystretch}{1.4}
\begin{array}{ll}
\text{Precision} &  80 \\
\text{Number of terms} &   200\\
\Gamma^\ast_0(4) & \big< T, \begin{smat}
{1}{-1}{4}{-3}
\end{smat} \big>\\
\om_1 &  i\cdot 0.000960627675025996\dots\\
\om_2 &   -0.02998129737318938\dots   \\
y_f    &    y_{-7}  \\
t      &   1 \\
D_1   & -7

\end{array}
\]

Similar to the last example, we chose the first positive residue modulo $2N$ which satisfies $D\equiv r^2\bmod 4N$ for each $D$.

\begin{table}[h]
\begin{center}
\begin{tabular}{|rr | rr|}\hline
\rule[-6pt]{0pt}{19pt}
$|D|$ & $m_{D,r}$ & $|D|$ & $m_{D,r}$ \\\hline
7 & 1 &103 & 1649 \\
15 & 5 & 111 & -765 \\
23 & -3 & 119 & -90 \\
31 & -50 &127 & 2664 \\
39 & -35 & 143 & -3729 \\
47 & 186 & 151 & -505 \\
55 & 215 & 159 & -2825 \\
71 & -315 & 167 & 3819 \\
79 & -10 & 183 & 2539 \\
87 & -497 & 191 & 1830 \\
95 & 405 & 199 & -5755 \\\hline
\end{tabular}
\vspace{.1in}
\caption{$f\in S_{12}(4)$. List of $D$, $m_{D,r}$ such that $ y_{D,r} - m_{D,r}y_f \in L$ for $|D|<200$.}\label{Ta:T5}
\end{center}
\end{table}

\end{ex}

%----------------------------------------------------------------

\section*{Acknowledgments}

I am deeply grateful to my advisor, Fernando Rodriguez Villegas, for his continuing guidance and support and for sharing his ideas that have enriched this work. I would also like to thank Don Zagier, Winfried Kohnen, and Chad Schoen for their ideas and suggestions. Thanks to Jeffrey Stopple for his careful reading of the manuscript.
The coefficients of modular forms were computed with the help of William Stein. This research was partially funded by the Donald D. Harrington Endowment Fellowship.

% ----------------------------------------------------------------

\bibliographystyle{akpbib}
\bibliography{mybibliography}

\end{document}